\documentclass{amsart}
\usepackage{amsmath}
\usepackage{amssymb}
\usepackage[all]{xy}

\newtheorem{theorem}{Theorem}[section]
\newtheorem{lemma}[theorem]{Lemma}
\newtheorem{prop}[theorem]{Proposition}

\newtheorem*{prop*}{Proposition}
\newtheorem*{lemma*}{Lemma}

\theoremstyle{remark}

\newtheorem{remarks}[theorem]{Remarks}

\newtheorem*{thmcomment}{Background}

\newtheorem*{thmcomment2}{Hooptedoodle}

\def\T{{\mathbb T}}

\def\Z{{\mathbb Z}}

\newcommand{\Aut}{\operatorname{Aut}}

\numberwithin{equation}{section}

\begin{document}
\SelectTips{cm}{}
\objectmargin={1pt}
\title{Contractible subgraphs and Morita equivalence of graph $C^*$-algebras}

\author{Tyrone Crisp}\author{Daniel Gow}

\address{School of Mathematical and Physical Sciences, The University of Newcastle, Callaghan, NSW 2308, Australia}
\email{tyrone.crisp@studentmail.newcastle.edu.au}

\address{School of Mathematics, The University of New South Wales, Sydney NSW 2052,
Australia} \email{danielg@maths.unsw.edu.au}

\date{February 2004}

\thanks{This research was supported by grants from the Australian Research
Council. We thank Iain Raeburn of the University of Newcastle for
helping us obtain this support.}


\maketitle

\begin{abstract}
In this paper we describe an operation on directed graphs which
produces a graph with fewer vertices, such that the $C^*$-algebra
of the new graph is Morita equivalent to that of the original
graph. We unify and generalize several related constructions,
notably delays and desingularizations of directed graphs.
\end{abstract}

\section{Introduction}

In recent years several authors have investigated certain
constructions on directed graphs, derived from the theory of
topological dynamics, which preserve the Morita equivalence class
of the associated graph $C^*$-algebras
\cite{Ash},\cite{B},\cite{BP},\cite{D},\cite{DT}. Typically, these
constructions have been viewed as enlargements of a graph which
preserve its path structure. In \cite{BP}, an attempt was made to
unify and generalize some of the existing results on the subject,
including the idea of a \emph{delay}, which was examined in
\cite{D} and is the basis of the desingularization of \cite{DT}.
It was noted \cite[Remarks 4.6]{BP} that applying a delay replaces
a vertex by a certain type of tree (called a \emph{gantlet} in
\cite{D}), and that the Morita equivalence results for delays may
still hold when a vertex is replaced by a more general tree.

In this paper we consider the reverse question: we aim to describe
sufficient conditions for a subgraph of an arbitrary directed
graph to be contractible, in the sense that its vertex set may be
reduced to yield a graph whose $C^*$-algebra is Morita equivalent
to that of the original graph. After a brief review of the
established definitions and notations for graph algebras, we state
and prove our main result, Theorem \ref{thm}, which shows that any
finite tree is contractible in the above sense, and that a similar
construction may be applied to more general acyclic subgraphs. In
particular, this theorem combines several of the separate results
of \cite{BP} and may cover some further examples as well. Our
proof follows closely the direct methods of \cite{BP}, and makes
use of a powerful theorem of \cite{BHRS}, the gauge-invariant
uniqueness theorem. Proposition \ref{G^0} gives equivalent
conditions to those of Theorem \ref{thm} which may make the
theorem easier to apply, and which give some idea of what the
contracted graph will look like. Section \ref{examples} discusses
the relationship of our theorem to the existing results, and
provides examples of graphs which are not contractible.

\section{Preliminaries} A \emph{directed graph} $E=(E^0,E^1,r,s)$
consists of countable sets $E^0$ (vertices) and $E^1$ (edges), and
maps $r,s:E^1\to E^0$ describing the range and source of each
edge. A vertex which emits no edges is called a \emph{sink}; a
vertex which emits infinitely-many edges is called an
\emph{infinite emitter}. A graph which contains no infinite
emitters is called \emph{row-finite}. Sinks and infinite emitters
are collectively described as \emph{singularities}, and we denote
by $E^0_{sing}$ the set of all singularities in $E^0$. A
\emph{Cuntz-Krieger $E$-family} consists of mutually orthogonal
projections $\{P_v:v\in E^0\}$ and partial isometries $\{S_e:e\in
E^1\}$ with mutually orthogonal ranges satisfying

(a) $S_e^*S_e=P_{r(e)}$,

(b) $S_eS_e^*\leq P_{s(e)}$, and

(c) $P_v=\sum_{s(e)=v} S_e S_e^*$ if $v$ is not a singularity.

\noindent The graph $C^*$-algebra of $E$, denoted $C^*(E)$, is the
universal $C^*$-algebra generated by a Cuntz-Krieger $E$-family
$\{s_e,p_v\}$. For any graph $C^*$-algebra there is an action
$\gamma:\T\to \Aut(C^*(E))$  characterized by $\gamma_z(p_v)=p_v$
and $\gamma_z(s_e)=zs_e$ for $v\in E^0$ and $e\in E^1$. This
\emph{gauge action} is equivalent to the universal property of
$C^*(E)$:

\begin{theorem}\cite[Theorem $2.1$]{BHRS} Let $E$ be a directed graph, $\{S_e,P_v\}$ be a Cuntz-Krieger $E$-family and
$\pi:C^*(E)\to C^*(S_e,P_v)$ the homomorphism satisfying
$\pi(s_e)=S_e$ and $\pi(p_v)=P_v$. Suppose that each $P_v$ is
non-zero, and that there is a strongly continuous $\T$-action
$\alpha$ on $C^*(S_e,P_v)$ such that
$\alpha_z\circ\pi=\pi\circ\gamma_z$ for all $z\in\T$. Then $\pi$
is faithful.\end{theorem}

Using the standard definitions and notations for paths in $E$ and
the convention that a vertex is a path of length zero, we denote
the set of finite paths by $E^*$, and the set of infinite paths by
$E^{\infty}$. A finite path $\alpha$ of positive length is a
\emph{cycle} if $s(\alpha)=r(\alpha)$ and $s(\alpha_i)\neq
s(\alpha_j)$ for $i\neq j$. An acyclic infinite path $\mu$ is a
\emph{tail} if each $s(\mu_i)$ emits only $\mu_i$ and each
$r(\mu_i)$ receives only $\mu_i$. For $u,v\in E^0$ we say that
$u\ge v$ if there is path in $E^*$ from $u$ to $v$. For
$U\subseteq E^0$, we say $U\geq v$ if there exists $u\in U$ such
that $u\geq v$. We define $v\geq U$ in a similar manner.

For $X\subset E^0$, we denote by $\Sigma H(X)$ the smallest
saturated hereditary subset of $E^0$ containing $X$, as defined in
\cite[Remark 3.1]{BHRS}: $\Sigma H(X):=\bigcup_{n\geq 0}
\Sigma_n(X)$, where $\Sigma_n(X)$ is defined inductively by
\begin{align*}\Sigma_0(X)&:=X\cup \{w\in E^0 : X\geq w\},\\
\Sigma_{n+1}(X)&:=\Sigma_n(X)\cup \{w\in E^0 :
0<|s^{-1}(w)|<\infty \mbox{ and } s(e)=w \mbox{ imply } r(e)\in
\Sigma_n(X)\}.
\end{align*}
\cite[Section 3]{BHRS} demonstrated a correspondence between the
saturated hereditary subsets of $E^0$ and the gauge-invariant
ideals of $C^*(E)$, and \cite[Lemma 2.2]{BP} gives a way to find
full projections in $\mathcal{M}(C^*(E))$ by examining $\Sigma
H(X)$ for suitable $X$.

\section{Contractible subgraphs}
\begin{theorem}\label{thm} Suppose $E$ is a directed graph with no tails,
and suppose $G^0\subset E^0$ such that $E^0_{sing}\subseteq G^0$ and the
subgraph $T$ of $E$ defined by $T^0:=E^0\setminus G^0$ and $T^1:=\{e\in
E^1:s(e),r(e)\in T^0\}$ is acyclic. Suppose that
\begin{itemize}\item[(a)] each vertex in $G^0$ is the source of at most one infinite path $\rho\in E^\infty$ such that $r(\rho_i)\in T^0$ for $i\geq1$;
\end{itemize}
 and that for each $\mu\in  T^\infty$, \begin{itemize}\item[(b)] $G^0\geq s(\mu)$; \item[(c)] $|r^{-1}(s(\mu_i))|=1$ for all $i$; and \item[(d)] $e\in E^1, r(e)=s(\mu)$ implies $|s^{-1}(s(e))|<\infty$.\end{itemize} Let $G$ be the graph with vertex set $G^0$ and one edge $e_\beta$ for each path $\beta\in E^*\setminus E^0$ with
$s(\beta),r(\beta)\in G^0$ and $r(\beta_i)\in T^0$ for $1\leq i<|\beta|$,
such that $s(e_\beta)=s(\beta)$ and $r(e_\beta)=r(\beta)$. Then $C^*(G)$ is isomorphic to a full corner of $C^*(E)$.
\end{theorem}

\begin{remarks}  {\it Conditions:} First note that any graph $C^*$-algebra
can be approximated by the $C^*$-algebra of a graph with no tails, by
replacing each tail with a sink as in \cite[Lemma $1.2$]{BPRS}. Now we
specify an acyclic subgraph $T$ of $E$ containing none of the
singularities of $E$, such that (a) holds, and every infinite path in $T$ satisfies
conditions ($b$)--($d$). Let $v$ be a vertex on such a path
$\mu$. Condition ($b$) says that there is a path from $G^0$ to $v$. Condition ($c$) says that $v$ receives exactly one edge. Condition ($d$) says that if $u$ is a vertex which emits an edge with range $v$, then $u$ emits only finitely many edges.
{\it Construction:} The idea is to replace each path $\beta$ through
$T$ from $G^0$ to $G^0$ by a single edge $e_\beta$ having the same source
and range as $\beta$. Since $T$ is acyclic and contains no singularities,
it is reasonable to expect that this construction preserve the ideal
structure of the graph algebra.
\end{remarks}

In order to prove Theorem \ref{thm} we shall first need to
establish some important properties of the subgraph $T$:

\begin{lemma}\label{T>G} Suppose $v\in T^0$. Then $v\geq G^0$.\end{lemma}
\begin{proof} Suppose that no path in $E^*$ with source $v$ has range in
$G^0$. Notice that we now have that $v\geq u\Longrightarrow u\ngeq
G^0$. Now $v$ is not a sink, so it emits an edge $\mu_1$. By our
assumption, $r(\mu_1)\in T^0$ and $r(\mu_1)\ngeq G^0$. Continuing
in this manner gives a path $\mu=\mu_1\mu_2\ldots\in T^\infty$. By
Condition ($b$) we can find a path $\alpha\in E^*$ such that
$s(\alpha)\in G^0$, $r(\alpha)=v$ and $r(\alpha_i)\in T^0$ for
$1\leq i\leq |\alpha|$. Now the path $\mu$ is not a tail, and each
$s(\mu_i)$ receives only one edge (by Condition ($c$)), so in
particular there exists $k$ such that $s(\mu_k)$ emits an edge
$\nu_1$ distinct from $\mu_k$. As before we must have $r(\nu_1)\in
T^0$ and $r(\nu_1)\ngeq G^0$, so we can construct an infinite path
$\nu=\nu_1\nu_2\ldots\in T^\infty$ such that $s(\nu)=s(\mu_k)$.
Now $\alpha\mu$ and $\alpha\mu_1\ldots\mu_{k-1}\nu$ are distinct
infinite paths which contradict Condition ($a$). \end{proof}

\begin{lemma}\label{B_v} For $v\in E^0$ define $B_v=\{\beta\in
E^*\setminus E^0: s(\beta)=v, r(\beta)\in G^0 \textrm{ and }
r(\beta_i)\in T^0\textrm{ for } 1\leq i<|\beta|\}$.  Then:
\begin{itemize}
\item[(a)] Suppose that $v\in E^0$ and that $\alpha,\beta\in B_v$. Then
neither of $\alpha$ and $\beta$ is a proper extension of the other.
\item[(b)] Suppose $\mu\in E^*$ and $s(\mu), r(\mu)\in G^0$. Then $\mu$
is a product of paths in $\bigcup_{v\in G^0} B_v$.
\item[(c)] $B_v=\emptyset \Longrightarrow v\in G^0$
\item[(d)] Suppose $v\in T^0$. If $|B_v|$ is infinite then there exists $\mu\in T^\infty$ such that
$v=s(\mu)$.
\end{itemize}
\end{lemma}

\begin{proof}(a) Suppose $\alpha=\beta\gamma$ for some $\gamma\in E^*\setminus
E^0$. We have $|\beta|<|\alpha|$ and $r(\alpha_{|\beta|})=r(\beta)\in
G^0$, contradicting $\alpha\in B_v$.

(b) Let $\mu=\mu_1\mu_2\dots\mu_n$. The proof is by induction on
$n$. The basis step $n=1$ is given by definition of $B_v$. Now
suppose that the assertion holds for paths of length less than
$n$. If there exists $k$ such that $1\leq k<n$ and $r(\mu_k)\in
G^0$, then $\mu=(\mu_1\dots\mu_k)(\mu_{k+1}\dots\mu_n)$ and hence
$\mu$ is a product of paths in $\bigcup_{v\in G^0}B_v$ by the
inductive hypothesis. If no such $k$ exists, then $\mu\in
B_{s(\mu)}$ by definition, and $s(\mu)\in G^0$.

(c) Suppose that $v\in T^0$. By Lemma \ref{T>G} there is a path
$\alpha \in E^*$ with $s(\alpha)=v, r(\alpha)\in G^0$. Let $k$ be the
first positive integer such that $r(\alpha_k)\in G^0$. Then
$\alpha_1 \alpha_2 \ldots \alpha_k \in B_v$ by definition of $B_v$, hence
$B_v$ is non-empty.

(d) Suppose $B_v$ is infinite. The vertex $v$ emits only finitely
many edges, since $v\in T^0$. By the pigeonhole principle, at least one
such edge $\mu_1$ must be the first edge in infinitely many $\beta$ in
$B_v$. If $r(\mu_1)\in G^0$, then by definition we would have $\mu_1\in
B_v$, and part (a) then implies that
$\mu_1$ can have no proper extension in $B_v$. Since we chose $\mu_1$
such that it had infinitely many extensions in $B_v$, we must have
$r(\mu_1)\in T^0$. Notice that by definition, $\mu_1\in T^1$. Since
$B_{r(\mu_1)}$ is infinite, we can find an edge
$\mu_2\in s^{-1}(r(\mu_1))\cap T^1$ such that $r(\mu_2)\in T^0$ and
$B_{r(\mu_2)}$ is infinite. Repeating this construction gives a sequence
$(\mu_1,\mu_2,\ldots)$ of edges in $T^1$ such that for all $i$,
$s(\mu_{i+1})=r(\mu_i)\in T^0$. Then $\mu:=\mu_1 \mu_2 \ldots$ is the
infinite path required.
\end{proof}

\begin{lemma}\label{B_v=sinversev} Suppose that $v\in E^0$ and
$\sup\{|\beta|:\beta\in B_v\}=1$. Then $B_v=s_E^{-1}(v)$.
\end{lemma}

\begin{proof} Each path in $B_v$ must be a single edge, so we
have  $B_v\subseteq s_E^{-1}(v)$. For the reverse inclusion, suppose
$e\in s_E^{-1}(v)$ and $r(e)\in T^0$. Lemma \ref{B_v} (c) implies that
$B_{r(e)}$ contains a path $\beta$, so $e\beta \in B_v$ is a path of
length at least $2$, a contradiction. Hence we must have $r(e)\in G^0$, and so $e\in B_v$
by definition. \end{proof}

\begin{lemma}\label{B_v2} Suppose that $v\in E^0$ and $B_v$ is finite and
non-empty. Then
$0<|s_E^{-1}(v)|<\infty$, and if $\{s_e,p_w\}$ is a Cuntz-Krieger
$E$-family, we have \begin{equation*}\label{B_vid} p_v=\sum_{\beta \in
B_v} s_\beta s_\beta^*\end{equation*}
\end{lemma}

\begin{proof} We begin by proving that $v$ is non-singular. Now
$E^0_{sing}\subseteq G^0$, so we need only consider $v\in G^0$. Trivially $v$
cannot be a sink, so it remains to show that $|s_E^{-1}(v)|<\infty$. We do
this by showing that each edge in $s_E^{-1}(v)$ has an extension in
$B_v$. Suppose $s(e)=v$. If $r(e)\in G^0$ then by definition $e\in B_v$,
so suppose $r(e)\in T^0$. There exists $\beta \in B_{r(e)}$ by Lemma
\ref{B_v} (c), so $e\beta \in B_v$ is an extension of $e$. Hence
$|s_E^{-1}(v)|\leq|B_v|<\infty$.

Now define $N(v):=\sup\{|\beta|:\beta\in B_v\}$, and notice that this
number is well-defined whenever $B_v$ is finite and non-empty. We shall
prove the equality $p_v=\sum_{\beta \in B_v} s_\beta s_\beta^*$ by
induction on $N(v)$.

Suppose $N(v)=1$. Lemma \ref{B_v=sinversev} implies $B_v=s^{-1}_E(v)$, and
the Cuntz-Krieger relation at $v$ gives
$p_v=\sum_{f\in s^{-1}_E(v)}s_f s_f^*=\sum_{\beta
\in B_v} s_\beta s_\beta^*$.

Now suppose that $N(v)=k$ and that $p_w=\sum_{\beta \in B_w} s_\beta
s_\beta^*$ for any $w\in E^0$ with $1\leq N(w)< k$. Then once again $v$
is non-singular, so
\begin{align*} p_v&=\sum_{f\in s^{-1}(v)} s_f s_f^*=\sum_{f\in
s^{-1}(v)\cap r^{-1}(G^0)}s_f s_f^* + \sum_{f\in s^{-1}(v)\cap
r^{-1}(T^0)}s_f s_f^*\\  &=\sum_{\beta\in B_v\cap E^1}s_\beta
s_\beta^*+\sum_{f\in s^{-1}(v)\cap r^{-1}(T^0)}s_f p_{r(f)} s_f^*.
\end{align*}

Consider an edge $f\in s^{-1}(v)\cap r^{-1}(T^0)$. We must have
$B_{r(f)}$ non-empty by Lemma \ref{B_v} (c). Each path $\alpha\in
B_{r(f)}$ gives a path $f\alpha\in B_v$, so $B_{r(f)}$ must be finite and
satisfy $N(r(f))\leq k-1$. Furthermore, every $\beta\in B_v$ with
$|\beta|\geq 2$ has the form $g\alpha$ for some $g\in s^{-1}(v)\cap
r^{-1}(T^0\setminus G^0),\alpha\in B_{r(f)}$ : by definition
$r(\beta_1)\in s^{-1}(v)\cap r^{-1}(T^0)$ and
$\beta_2\ldots\beta_{|\beta|}\in B_{r(\beta_1)}$. Thus for each such $f$
we can apply the inductive hypothesis to $r(f)$, giving
\begin{align*} p_v&=\sum_{\beta\in B_v\cap E^1}s_\beta
s_\beta^*+\sum_{f\in s^{-1}(v)\cap r^{-1}(T^0)}\left(\sum_{\alpha\in
B_{r(f)}}s_{f\alpha} s_{f\alpha}^*\right)\\  &=\sum_{\beta\in B_v\cap
E^1}s_\beta s_\beta^*+\sum_{\substack{\beta\in B_v,\\|\beta|\geq
2}}s_\beta s_\beta^*\\  &=\sum_{\beta\in B_v} s_\beta s_\beta^*.
\end{align*} This completes the proof by induction.
\end{proof}

\begin{proof}[Proof of Theorem \ref{thm}] Let $\{s_e,p_v\}$ be the
canonical Cuntz-Krieger $E$-family that generates $C^*(E)$. For
$e_\beta\in G^1$ we define $T_{e_\beta}=s_\beta$, and for $v\in G^0$ we
define $Q_v=p_v$.

The $Q_v$ are mutually orthogonal projections because the $p_v$
are. The $T_{e_\beta}$ are partial isometries because they are
products of the partial isometries $s_f$ (recall the properties of
Cuntz-Krieger families of partial isometries). To see that they
have mutually orthogonal ranges, suppose $e_\alpha,e_\beta\in
G^1$, $e_\alpha\neq e_\beta$. Then $\alpha,\beta\in\bigcup_{v\in
G^0}B_v$ have the property that neither one is an extension of the
other, by Lemma \ref{B_v}(a). \cite[Lemma 1.1]{KPR}, which applies
to infinite graphs, then implies $s^*_\alpha s_\beta=0$, so
$T_{e_\alpha}T_{e_\alpha}^*T_{e_\beta}T_{e_\beta}^*=s_\alpha
s_\alpha^* s_\beta s_\beta^*=0$, and $T_{e_\alpha}$ and
$T_{e_\beta}$ have mutually orthogonal ranges. For $e_\beta\in
G^1$ we have $T^*_{e_\beta}
T_{e_\beta}=s_\beta^*s_\beta=p_{r(\beta)}=Q_{r(e_\beta)}$ and
$T_{e_\beta} T^*_{e_\beta}=s_\beta s_\beta^*\leq
p_{s(\beta)}=Q_{s(e_\beta)}$. Now suppose $v\in G^0$ is
non-singular in $G$: that is, suppose $0<|s_G^{-1}(v)|<\infty$.
Since $s_G^{-1}(v)$ is equinumerous with $B_v$, $B_v$ is finite
and non-empty. Lemma \ref{B_v2} then gives
\[Q_v=p_v=\sum_{\beta\in B_v}s_\beta s^*_\beta=\sum_{e_\beta\in
s_G^{-1}(v)}T_{e_\beta}T_{e_\beta}^*.\]  Thus $\{T_{e_\beta},Q_v\}$ is a
Cuntz-Krieger $G$-family.

A slightly modified form of the argument of \cite[Section 2]{BPRS} shows
that there is a strongly continuous action $\alpha$ of $\T$ on $C^*(E)$
such that $\alpha_z(p_v)=p_v$ for all $v\in E^0$ and
\begin{equation*}
\alpha_z(s_e)=\begin{cases} zs_e \textrm{ if $r(e)\in G^0$}\\ s_e
\textrm{ if $r(e)\in T^0$}\end{cases}
\end{equation*}

Let $\{t_{e_\beta},q_v\}$ be the canonical Cuntz-Krieger $G$-family. The
universal property of $C^*(G)$ ensures the existence of a homomorphism
$\pi$ of $C^*(G)$ onto $C^*(T_{e_\alpha},Q_v)$ such that
$\pi(t_{e_\beta})=T_{e_\beta}$ and $\pi(q_v)=Q_v$ for all $e_\beta\in
G^1$ and $v\in G^0$. If $\gamma$ is the canonical gauge action on
$C^*(G)$, then $\pi\circ\gamma_z(q_v)=\alpha_z\circ\pi(q_v)$ for all
$v\in G^0$ and $z\in\T$. Now fix $z\in\T$ and suppose $e_\beta\in G^1$.
Then $\pi\circ\gamma_z(t_{e_\beta})=\pi(zt_{e_\beta})=zT_{e_\beta}$, and
by definition of $B_v$,
$\alpha_z\circ\pi(t_{e_\beta})=\alpha_z(s_{\beta_1}\dots
s_{\beta_{|\beta|}})=s_{\beta_1}\ldots
s_{\beta_{|\beta|-1}}zs_{\beta_{|\beta|}}=zT_{e_{\beta}}$. Hence
$\pi\circ\gamma=\alpha\circ\pi$ on all of $C^*(G)$. The gauge-invariant
uniqueness theorem \cite[Theorem 2.1]{BHRS} then implies that $\pi$ is an
isomorphism of $C^*(G)$ onto $C^*(T_{e_\beta},Q_v)$. We prove Theorem
\ref{thm} by showing that $C^*(T_{e_\beta},Q_v)$ is a full corner of
$C^*(E)$.

By \cite[Lemma 1.2(c)]{BPRS} the sum $\sum_{v\in G^0} p_v$ converges
strictly to a projection $P\in \mathcal M (C^*(E))$. We claim that
$C^*(T_{e_\beta},Q_v)=PC^*(E)P$.

Note that for each $v\in G^0$ we have $Q_v\leq P$, so $Q_v=PQ_v P\in
PC^*(E)P$. We then have that for every $e_\beta\in G^1$,
\[T_{e_\beta}=Q_{s(e_\beta)}T_{e_\beta}Q_{r(e_\beta)}=PQ_{s(e_\beta)}T_{e_\beta}Q_{r(e_\beta)}P\in
PC^*(E)P.\] So $C^*(T_{e_\beta},Q_v)\subseteq PC^*(E)P$ is easy.

Now fix $s_\mu s_\nu^*\in C^*(E)$. Then $Ps_\mu s_\nu ^*P=\sum_{v,w\in
G^0} p_v s_\mu s_\nu ^*p_w$. Now $p_v s_\mu=0$ unless $v=s(\mu)$, in
which case $p_v s_\mu=s_\mu$. We can apply the same argument to
$s_\nu^*p_w=(p_w s_\nu)^*$, so we have the following: Suppose that
$Ps_\mu s_\nu ^*P\neq0$. Then $Ps_\mu s_\nu^*P=s_\mu s_\nu^*$, $s(\mu),
s(\nu)\in G^0$, and $r(\mu)=r(\nu)$. Hence to show $Ps_\mu s_\nu ^*P\in
C^*(T_{e_\beta},Q_v)$ it will suffice to consider the following three
cases:

\begin{enumerate}
\item $r(\mu)\in G^0$;
\item $r(\mu)\in T^0$ and $B_{r(\mu)}$ is finite; and
\item $r(\mu)\in T^0$ and $B_{r(\mu)}$ is infinite,
\end{enumerate} and to show in each case that $s_\mu s_\nu^*\in
C^*(T_{e_\beta},Q_v)$.

{\bf Case 1.} By Lemma \ref{B_v}(b) we can write $\mu$ as a product
$\alpha^1 \alpha^2 \ldots \alpha^n$ of paths in $\bigcup_{v\in G^0}B_v$,
so that $s_\mu=s_{\alpha^1}s_{\alpha^2}\ldots s_{\alpha^n}\in
C^*(T_{e_\beta},Q_v)$. Similarly we can write $s_\nu =
s_{\beta^1}s_{\beta^2}\ldots s_{\beta^m}\in C^*(T_{e_\beta},Q_v)$, so
$s_\mu s_\nu^* \in C^*(T_{e_\beta},Q_v)$.

{\bf Case 2.} First notice that Lemma \ref{B_v}(c) implies that
$B_{r(\mu)}$ is non-empty. Let $k:=\max\{i:s(\mu_i)\in G^0\}$, and
consider the paths $\rho:=\mu_1\ldots\mu_{k-1}$ and
$\gamma:=\mu_k\ldots\mu_{|\mu|}$. We can decompose $\mu$ as the product
$\rho\gamma$ such that
$s(\rho),r(\rho)\in G^0$, $s(\gamma)\in G^0$ and $r(\gamma_i)\in T^0$ for
$1\leq i\leq |\gamma|$. Similarly, we may write
$\nu=\sigma\delta$ for some paths $\sigma,\delta\in E^*$ with the same
properties as $\rho$ and $\gamma$, respectively. Case $1$ shows that
$s_\rho$ and $s_\sigma$ are in $C^*(T_{e_\beta},Q_v)$, so to show $s_\mu
s_\nu^*\in C^*(T_{e_\beta},Q_v)$ it will be enough to show that $s_\gamma
s_\delta^*\in C^*(T_{e_\beta},Q_v)$. Since $B_{r(\mu)}$ is finite and
non-empty, and $r(\mu)=r(\gamma)=r(\delta)$, we can use Lemma \ref{B_v2}
to get
\begin{align*} s_\gamma s_\delta^*&=s_\gamma p_{r(\mu)} s_\delta^*\\  &=
s_\gamma \left(\sum_{\beta \in B_{r(\mu)}} s_\beta s_\beta^*
\right)s_\delta^*\\  &=\sum_{\beta \in B_{r(\mu)}} s_{\gamma\beta}
s_{\delta\beta}^* \\ &=\sum_{\beta \in B_{r(\mu)}} T_{e_{\gamma\beta}}
T_{e_{\delta\beta}}^*
\quad\textrm{ (since each $\gamma\beta, \delta\beta\in \bigcup_{v\in
G^0}B_v$)}.
\end{align*} Thus $s_\mu s_\nu^*$ is a finite sum of elements of
$C^*(T_{e_\beta},Q_v)$, so $s_\mu s_\nu^*\in   C^*(T_{e_\beta},Q_v)$.

{\bf Case 3.} (This is a combination of our method for Case 2 and the
proof of \cite[Theorem $4.2$]{BP}.) Lemma \ref{B_v}(d) implies that
$r(\mu)$ is the source of some infinite path $\epsilon\in T^\infty$. By
Condition ($b$) there exists a path $\alpha=\alpha_1\ldots\alpha_n\in E^*$
such that $s(\alpha)\in G^0$,
$\alpha_2\alpha_3\ldots\alpha_{|\alpha|}\in T^*$ and
$r(\alpha)=r(\mu)$. For convenience we shall write $\alpha_0$ to denote the
vertex $s(\alpha)$ viewed as a path of zero length, and for $0\leq i\leq
n$ we shall write $\gamma_i:=\alpha_0\ldots\alpha_i$. Condition ($c$)
implies that for each $i\geq 1$,
$r(\alpha_i)$ receives only one edge, so we must have $\mu=\rho\alpha$ and
$\nu=\sigma\alpha$ for some $\rho,\sigma\in E^*$. Once again Case $1$
shows that $s_\rho,s_\sigma\in C^*(T_{e_\beta},Q_v)$, so we need only
prove $s_\alpha s_\alpha^*\in C^*(T_{e_\beta}, Q_v)$.

By our definition, $s_{\alpha}s_{\alpha}^*=s_{\gamma_n}s_{\gamma_n}^*$;
our plan for proving $s_\alpha s_\alpha^*\in C^*(T_{e_\beta}, Q_v)$ is to
reduce this product to $p_{s(\alpha)}- A$, where $A$ is a finite sum of
elements of $C^*(T_{e_\beta}, Q_v)$. If we can do this, we will be done
with Case $3$: $s(\alpha)\in G^0$ implies
$p_{s(\alpha)}-A=Q_{s(\alpha)}-A\in C^*(T_{e_\beta}, Q_v)$. We perform
this reduction recursively, as follows:

Suppose $0\leq k<n$. We shall show that
$s_{\gamma_{k+1}}s_{\gamma_{k+1}}^*=s_{\gamma_k}s_{\gamma_k}^*-A_{k+1}$
where $A_{k+1}$ is a finite sum of elements of $C^*(T_{e_\beta}, Q_v)$.
Since $\alpha_2\ldots\alpha_n\epsilon$ is a path in $T^\infty$, Condition
($d$) implies that each $s^{-1}(s(\alpha_i))$ is finite. Suppose
$s^{-1}(s(\alpha_{k+1}))=\{\alpha_{k+1},f_1,\ldots,f_m\}$. Then
$s_{\alpha_{k+1}}s_{\alpha_{k+1}}^*=p_{s(\alpha_{k+1})}-\sum_{i=1}^m
s_{f_i}s_{f_i}^*$, and so \[s_{\gamma_{k+1}}s_{\gamma_{k+1}}^*
=s_{\gamma_k}s_{\alpha_{k+1}}s_{\alpha_{k+1}}^*s_{\gamma_k}^*
=s_{\gamma_k}p_{s(\alpha_{k+1})}s_{\gamma_k}^*-\sum_{i=1}^m
(s_{\gamma_k}s_{f_i} s_{f_i}^*s_{\gamma_k}^*).\] Since
$s(\alpha_{k+1})=r(\gamma_k)$,
$s_{\gamma_k}p_{s(\alpha_{k+1})}s_{\gamma_k}^*=s_{\gamma_k}s_{\gamma_k}^*$
so we would like to show that each $s_{\gamma_{k}}s_{f_i}
s_{f_i}^*s_{\gamma_{k}}^*$ is equal to a finite sum of elements of
$C^*(T_{e_\beta}, Q_v)$. Fix $f\in
s^{-1}(s(\alpha_{k+1}))\setminus\{\alpha_{k+1}\}$. If $r(f)\in G^0$, then
the path $\gamma_{k}f$ is in $B_{s(\alpha)}$, and hence
$s_{\gamma_{k}}s_{f_i}
s_{f_i}^*s_{\gamma_{k}}^*=T_{e_{\gamma_{k}f}}T_{e_{\gamma_{k}f}}^*\in
C^*(T_{e_\beta}, Q_v)$ as required. Now suppose $r(f)\notin G^0$. Then
$B_{r(f)}$ is non-empty by Lemma \ref{B_v}(c); suppose it is infinite.
Lemma \ref{B_v}(d) implies that $r(f)=s(\zeta)$ for some
$\zeta\in T^\infty$. Then $\gamma_k f\zeta$ and $\alpha \epsilon$ are distinct paths which contradict Condition ($a$). So
$B_{r(f)}$ must in fact be finite. Hence we can apply Lemma \ref{B_v2} to
give
\begin{align*} s_{\gamma_{k}}s_{f}
s_{f}^*s_{\gamma_{k}}^*&=s_{\gamma_{k}}s_{f}p_{r(f)}
s_{f}^*s_{\gamma_{k}}^*\\  &=\sum_{\beta\in
B_{r(f)}}s_{\gamma_{k}}s_{f}s_{\beta}s_\beta^*s_{f}^*s_{\gamma_{k}}^*.\\
\end{align*} Each path $\gamma_{k}f\beta$ satisfies the requirements for
membership of
$B_{s(\alpha)}$, so we have \[s_{\gamma_{k}}s_{f}
s_{f}^*s_{\gamma_{k}}^*=\sum_{\beta\in
B_{r(f)}}T_{e_{\gamma_{k}f\beta}}T_{e_{\gamma_{k}f\beta}}^*.\]  This is a
finite sum of elements of $C^*(T_{e_\beta}, Q_v)$, which is precisely
what we wanted.

We can now apply this reduction $n$ times to get $s_\alpha
s_\alpha^*=s_{\gamma_0}s_{\gamma_0}^*-\sum_{k=1}^n A_k$, and since
$s_{\gamma_0}=s_{\alpha_0}=p_{s(\alpha)}$ we are done with Case $3$.
\medskip

Now suppose $a=PAP\in PC^*(E)P$. We can find $A_n\in C^*(E)$ such that each
$A_n$ is a finite linear combination of elements of the form $s_\mu
s_\nu^*$ and $A_n\rightarrow A$. Cases $1$, $2$ and $3$ show that each
$PA_nP\in C^*(T_{e_\beta}, Q_v)$, and continuity of multiplication in
$\mathcal{M}(C^*(E))$ implies $PA_nP\rightarrow PAP=a$. This gives
$PC^*(E)P\subseteq C^*(T_{e_\beta}, Q_v)$, and it remains to show that
$P$ is a full projection.

By \cite[Lemma 2.2]{BP}, to show that the projection $P$ is
full we have only to show that $E^0\subset\Sigma H(G^0)$. We have
$G^0\subset\Sigma H(G^0)$ by definition, so suppose $v\in T^0$, and note
that Lemma \ref{B_v}(c) implies that $B_v$ is non-empty. If
$B_v$ is infinite, then $v=s(\mu)$ for some $\mu\in T^\infty$ by Lemma
\ref{B_v}(d); Condition ($b$) then implies $G^0\geq v$, so $v\in \Sigma
H(G^0)$. Now suppose $B_v$ is finite, and note that Lemma \ref{B_v2}
implies $0<|s^{-1}(v)|<\infty$. Define $N(v):=\sup\{|\beta|:\beta\in
B_v\}$ as in the proof of Lemma \ref{B_v2}. We show that $v\in \Sigma
H(G^0)$ by induction on $N(v)$.

Suppose $N(v)=1$. Lemma \ref{B_v=sinversev} shows that every edge in
$s^{-1}(v)$ has range in $G^0$, and since $0<|s^{-1}(v)|<\infty$ we have
$v\in \Sigma H(G^0)$. Now assume that for $w\in T^0$, $1\leq N(w)\leq k $
implies that $w\in \Sigma H(G^0)$, and suppose $N(v)=k+1>1$. For each
edge $e\in s^{-1}(v)$, we must have either $r(e)\in G^0$ or
$r(e)\in T^0$. Suppose $r(e)\in T^0$. Lemma \ref{B_v}(c) implies that $B_{r(e)}\neq\emptyset$, so $N(r(e))\geq 1$. Now any path $\beta\in B_{r(e)}$ gives a path $e\beta\in B_v$ with length $|\beta|+1$, so we must have $1\leq N(r(e))\leq k$. The inductive hypothesis then implies
$r(e)\in \Sigma H(G^0)$. Thus any edge $e$ with source $v$ has range in
$\Sigma H(G^0)$, and since $v$ is non-singular, $v$ is then in $\Sigma
H(G^0)$. It follows, by induction, that $T^0\subset\Sigma H(G^0)$. Hence
$E^0 \subset \Sigma H(G^0)$ as required, $P$ is a full projection in
$\mathcal{M}(C^*(E))$, and $C^*(G)$ is isomorphic to a full corner of
$C^*(E)$.
\end{proof}

The conditions of Theorem \ref{thm} are based on a description of
$T$, the subgraph to be contracted. It is possible to formulate
equivalent conditions based on a description of $G^0$, the vertex
set of the graph obtained by the contraction:

\begin{prop} \label{G^0}Suppose $E$ is a directed graph with no tails, and
suppose $G^0\subset E^0$ such that $E^0_{sing}\subseteq G^0$. Then
$G^0$ satisfies the conditions of Theorem \ref{thm} if and only if
it satisfies the following:
\begin{enumerate}
\item $\lambda\in E^*$ a cycle $\Longrightarrow s(\lambda_i)\in G^0$ for
some $i$; and
\item Suppose $\mu,\nu\in E^\infty$ are distinct and acyclic; then
\begin{itemize}

\item[($a'$)] $s(\mu)=s(\nu)\in G^0\Longrightarrow r(\mu_i)\in G^0$ or
$r(\nu_i)\in G^0$ for some $i$.
\item[($b'$)] Either $G^0\geq s(\mu)$ or $s(\mu_i)\in G^0$ for some $i$;
\item[($c'$)] $|r^{-1}(s(\mu))|>1\Longrightarrow s(\mu_i)\in G^0$ for some
$i$;
\item[($d'$)] $|s^{-1}(s(\mu))|=\infty\implies r(\mu_i)\in G^0$ for some
$i$.

\end{itemize}
\end{enumerate}
\end{prop}
\begin{proof} Suppose $\lambda\in E^*$ is a cycle. By definition of $T^1$,
$\lambda\in T^*\iff s(\lambda_i)\in T^0$ for all $i$. Now any
cycle in $T^*$ is a cycle in $E^*$, so the subgraph $T$ is acyclic
if and only if ($1$) holds. Now each of Conditions ($a$)--($d$) in
Theorem \ref{thm} is equivalent to the corresponding condition in
the proposition. For example, suppose ($b$) holds, and let $\mu\in
E^\infty$ with $s(\mu_i)\in T^0$ for all $i$. Then by definition
we have $\mu\in T^\infty$, so ($b$) implies $G^0\geq s(\mu)$,
giving ($b'$). Conversely, suppose ($b'$) holds, and let $\mu\in
T^\infty$. Then again by definition of $T^1$ we cannot have any
$s(\mu_i)\in G^0$, so $G^0\geq s(\mu)$, giving ($b$).
\end{proof}

\section{Examples}\label{examples} (i) For a  graph $E$ with no tails, the
desingularization $F$ of $E$, described in \cite{DT}, is obtained
by adding a tail at each infinite-emitter $v_0\in E^0$ and
distributing the edges in $s^{-1}(v_0)$ along this tail, such that
the resulting graph is row-finite. It is straightforward to check
that each such tail is acyclic, non-singular and satisfies
Conditions ($a$)--($d$) of Theorem \ref{thm}, and hence with
$G^0=E^0$, Theorem \ref{thm} gives Morita equivalence of $C^*(E)$
with $C^*(F)$ as in \cite[Theorem $2.11$]{DT}.

(ii) \cite[Section 5]{B} examined a relation on directed graphs
called \emph{elementary strong shift equivalence}. Two graphs
$E_1$ and $E_2$ are elementary strong shift equivalent via $E_3$
if $E_3$ is a bipartite graph whose vertex set $E_3^0$ is the
disjoint union $E_1^0\cup E_2^0$ such that the paths of length $2$
in $E_3^*$ with source and range in $E_i^0$ are in one-to-one
correspondence with the edges in $E_i$. It was shown
(\cite[Theorem $5.2$]{B}) that for row-finite graphs with no
sinks, elementary strong shift equivalence implies Morita
equivalence of the associated graph algebras. Setting $E=E_3$ and
$G^0=E_i^0$ for $i=1,2$ in Theorem \ref{thm} gives the same
result. (To see that this is an applicable choice of $G^0$, notice
that $E^0_{sing}=\emptyset$ and that the corresponding subgraph
$T$ has no edges.)

(iii) An out-delay $d_s(E)$ of a graph $E$ as described in
\cite{BP} is obtained by adding some subpath of a tail (called a
\emph{gantlet} in \cite{D}), possibly of length zero or $\infty$,
to each vertex $v_0\in E^0$ and distributing the edges in
$s^{-1}(v_0)$ along this path. It can be seen that an out-delay is
strictly proper, as defined in \cite{BP}, if and only if the graph
$T$ defined as the union of all the added gantlets satisfies
Conditions ($a$)--($d$) of Theorem \ref{thm}. Now taking
$G^0=E^0$, our Theorem \ref{thm} gives Theorem $4.2$ of \cite{BP}.
In a similar way Theorem \ref{thm}, applied to an in-delay
$d_r(E)$ of a graph $E$, gives \cite[Theorem 4.5]{BP} when
$G^0=E^0$. The equivalence theorem for in-splittings (Corollary
5.4 of \cite{BP}) now follows from \cite[Theorem 5.3]{BP} and
Theorem \ref{thm}.

(iv) Theorem \ref{thm} may be more general than the results of
\cite{BP}; in particular, it covers situations where it may not be
obvious that a finite sequence of delays and splittings give the required
reduction. For example, denote by $B_n$ the binary tree with $n$ generations and all edges directed toward the leaves. Now let $\xygraph{ !~:{@{-}|@{>}}v:[r]{B_n}:[r]w}$ denote the graph with one edge from $v$ to the root of $B_n$, and one edge from each leaf of $B_n$ to $w$. For example, for $n=2$ the graph is
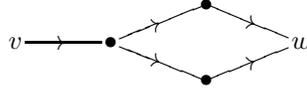
\[\xygraph{ !~:{@{-}|@{>}}
{v}:[r]="root"\bullet([r]([u(0.4)]="top"\bullet,[d(0.4)]="bottom"\bullet,[r]{w}="w"),:"top":"w",:"bottom":"w")}\]
Now consider the following graph $E$:
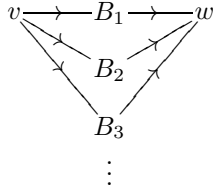
\[\xygraph{
!~:{@{-}|@{>}}
v(:[r]{B_1}([d(0.6)]{B_2}[d(0.6)]{B_3}[d(0.4)]\vdots,:[r]w),:"B_2":"w",:"B_3":"w"}\]
Letting $G^0:=\{v,w\}$, it can be seen that the conditions of
Theorem \ref{thm} are satisfied and hence that $C^*(E)$ is Morita
equivalent to $C^*(\xygraph{!~:{@{-}|@{>}}v:^{\infty}[r]w})$. It
is not obvious that this result could be deduced from
finitely-many applications of the results in \cite{BP}. Indeed, it
seems reasonable to assume that the smallest number $r_n$ of
applications of those results required to deduce the equivalence
of $C^*\left(\xygraph{ !~:{@{-}|@{>}}v:[r]{B_n}:[r]w}\right)$ and
$C^*(\xygraph{!~:{@{-}|@{>}} v:^{2^{n-1}}[r]w})$ should increase
without bound as $n$ tends to infinity, and that the number of
such applications required to deduce equivalence of $C^*(E)$ and
$C^*(\xygraph{!~:{@{-}|@{>}}v:^{\infty}[r]w})$ should exceed every
$r_n$.

(v) In \cite{HS2} the $C^*$-algebra $C(L_q(p;m_1,\dots,m_n))$ of
continuous functions on the quantum lens space
$L_q(p;m_1,\dots,m_n)$ was defined as the fixed point algebra
$C(S_q^{2n-1})^{\tilde{\Lambda}}$ of a certain action of $\Z_p$ on
the $C^*$-algebra of continuous functions on the odd-dimensional
quantum sphere $S_q^{2n-1}$. It was shown in \cite[Theorem
4.4]{HS1} that $C(S_q^{2n-1})$ is isomorphic to the $C^*$-algebra
of a directed graph $L_{2n-1}$. If $\Lambda$ is the $\Z_p$-action
on this graph algebra corresponding to $\tilde{\Lambda}$, then
\cite[Corollary 2.5]{KP} shows that the crossed product
$C^*(L_{2n-1})\times_\Lambda \Z_p$ is itself the $C^*$-algebra of
a certain graph, called the \emph{skew product} graph. In
\cite[Theorem 2.5]{HS2}, the fixed point algebra corresponding to
$\Lambda$ was also realized as a graph algebra
$C^*\left(L_{2n-1}^{(p;m_1,\dots,m_n)}\right)$. It can be seen
that the graph $L_{2n-1}^{(p;m_1,\dots,m_n)}$ may be obtained from
the skew product graph $L_{2n-1}\times_c \Z_p$ by a contraction as
in Theorem \ref{thm} (indeed, take $G^0=(L_{2n-1})^0\times
\{0\}$). It is likely that this result can be generalized, so that
the fixed point algebras corresponding to certain actions of
finite groups on graph algebras may themselves be realized as
graph algebras, with the graph in question being obtained from the
skew product graph by a contraction as in this example.

(vi) The following examples illustrate cases where Theorem
\ref{thm} is not applicable. First, consider the following graph
$E$:

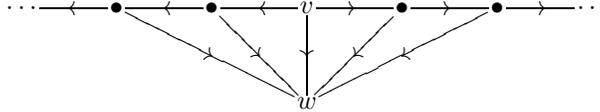
\[\xygraph{ !~:{@{-}|@{>}} v(:[l]\bullet="l1":[l]\bullet="l2":[l]
\cdots,:[r]\bullet="r1":[r]\bullet="r2":[r]\cdots,:[d]w)"l1":"w"
"l2":"w" "r1":"w" "r2":"w" }\] If we relax Condition ($a$), we can
take $G^0=\{v,w\}$ and deduce Morita equivalence of $C^*(E)$ and
$C^*(\xygraph{v:@{-}|@{>}^{\infty}[r]w})$. However, using the
results of \cite{BHRS} it can be seen that $C^*(E)$ has $3$
non-trivial ideals, while
$C^*(\xygraph{v:@{-}|@{>}^{\infty}[r]w})$ has only one. Now
consider the following graph $F$:

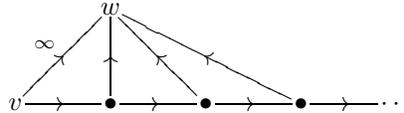
\[\xygraph{!~:{@{-}|@{>}} v:[r]\bullet(:[u]w):[r]\bullet(:"w"):[r]
\bullet(:"w"):[r]\cdots "v":^\infty "w"}\] Relaxing Condition
($d$) and taking $G^0=\{v,w\}$ gives Morita equivalence of
$C^*(F)$ and $C^*(\xygraph{v:@{-}|@{>}^{\infty}[r]w})$, which is
again contradicted by counting the saturated hereditary subsets of
$F^0$.

\end{document}